\documentclass[a4paper, 12pt]{article}

\usepackage[koi8-r,cp1251]{inputenc}
\usepackage{amsfonts,amssymb,amsmath, hyperref}
\usepackage[final]{epsfig}
\usepackage{graphicx}

\begin{document}

\begin{center}
{\Large On the boundedness of  Hausdorff  operators
on real Hardy spaces $H^1$ over homogeneous  spaces of groups with local doubling property}
\end{center}

\begin{center}
A. R. Mirotin
\end{center}

\begin{center}
amirotin@yandex.ru
\end{center}

\textsc{Abstract}. We give conditions for boundedness of
Hausdorff operators on real Hardy spaces $H^1$ over homogeneous spaces of
locally compact groups with local doubling property. The special  case of the hyperbolic plane is considered.

Key wards. Hausdorff operator, locally compact group, homogeneous
space, atomic Hardy space,  hyperbolic plane.

MSC classes:	47G10, 43A85, 47A30, 51M10, 22E30

\section{Introduction}

Hausdorff  operators on the finite interval were introduced by Hardy \cite[Chapter XI]{H}.
Its natural $n$-dimensional generalization due to Lerner and Liflyand \cite{LL}. For more information on the history of the issue, see the survey   articles  \cite{Ls}, \cite{CFW}.

Below we shall denote by $\mathrm{Aut}(G)$  the space  of all topological  automorphisms of a topological group $G$ endowed with its natural topology (see, e.g. \cite[(26.3)]{HiR}).

In \cite{JMAA}, \cite{AddJMAA} the next definition was proposed \footnote{The  case of a Hausdorff  operator on  $p$-Adic vector spaces  was considered earlier in \cite{Vol}, the special case of a Hausdorff  operator on the Heisenberg group in the sense of this definition was considered   in \cite{RFW}, and \cite{CFWu}.}

 \textbf{Definition}. Let  $(\Omega,\mu)$ be  a measure space, $G$ a topological group, $A:\Omega\to \mathrm{Aut}(G)$ a measurable map,
   and $\Phi$ a locally $\mu$-integrable function on $\Omega.$
   We define the \textit{Hausdorff  operator} with the kernel $\Phi$ over the group $G$  by the formula
$$
(\mathcal{H}f)(x)=\int_\Omega \Phi(u)f(A(u)(x))d\mu(u).
$$

By \cite[Lemma 1]{JMAA} for a locally compact group  $G$ this operator is bounded on $L^p(G)$ ($1\le p\le\infty$) provided $\Phi(u)(\mathrm{mod} A(u))^{-1/p}\in L^1(\Omega,\mu)$ and
$$
\|\mathcal{H}\|_{L^p\to L^p}\le \int_\Omega |\Phi(u)|(\mathrm{mod} A(u))^{-1/p}d\mu(u).
$$
In  \cite[Theorem]{AddJMAA} sufficient conditions  were  given for boundedness of a Hausdorff  operator on atomic Hardy space $H^1(G)$.

In the case $G=\mathbb{R}$ there are many classical operators in analysis which are special cases of the Hausdorff operator  in the sense of previous definition for suitable measure spaces $(\Omega,\mu)$ (see, e.g., \cite{WF}, \cite{JJS}, \cite{S},  \cite{CDFZ} and the bibliography therein). The next  example shows that the  Harish-Chandra transform leads to a Hausdorff operator on the group $SO(2,\mathbb{R})$.

\textbf{Example 1}.
One of the form of the Harish-Chandra transform for the group  $G=SL(2,\mathbb{R})$  looks as follows: $H^Kf(x)=D(x)(Hf)(x)$\ ($x\in SO(2,\mathbb{R})$) where $D(x)=2i\sin\theta$ if the matrix $x$ represents the rotation of the plane by an angle $\theta$. Here $H$ stands for a Hausdorff  operator
$$
(Hf)(x)=\int_{A^+}\frac{\alpha(u)+\alpha(u^{-1})}{2}f(u^{-1}xu)d\mu(u),
$$
 where $A^+$ denotes the set of $2\times 2$ matrices of the form $u=\mathrm{diag}(a,a^{-1})$, $a\ge 1$, $\alpha(u)=a^2$, and $d\mu(u)=da$  (see, e.g., \cite[Chapter VII, \S 5]{Leng}).

In \cite{JMAA}, \cite{AddJMAA} classical results  on the boundedness of Hausdorff operators on the Hardy space $H^1$ over finite-dimensional real space were
generalized to the case of a Hardy space over locally compact metrizable groups
with the doubling property. In \cite{JBSU}  Hausdorff operators on Lebesgue and real Hardy spaces over
homogeneous spaces of locally compact groups with doubling property were
considered.  On the other hand, as was shown by T.~ Kawazoe,
 non compact semisimple Lie groups should not satisfy the doubling property
 but often enjoy the less restrictive so called local doubling property (see condition (LDP)
below)  \cite[Lemma 2.6]{Kaw},    \cite[Lemma 3.2]{KawL} \footnote{The author thanks Professor R. Daher for these references.}.
The aim of this work, is to give conditions for boundedness of
Hausdorff operators on real Hardy spaces $H^1$ over homogeneous spaces of
locally compact groups with local doubling property (in particular, over Riemannian symmetric spaces $G/K$, where $G$ is a connected semisimple Lie group
with finite center and $K$ its maximal compact subgroup).

\section{ The main result}

 We shall assume in this section that  $G$ is a locally compact metrizable group
 with left invariant metric $\rho$ and left Haar measure $\nu$, $K$ its  compact subgroup with normalized Haar measure $\beta$.

In \cite{JBSU} Hausdorff operators on the homogeneous space $G/K$ were introduced in the following way. Recall that the quotient space $G/K$ consists of left cosets $\dot x := xK = \pi_K(x)$ ($x\in G$)
where $\pi_K : G \to G/K$ stands for a natural projection.
 We assume that
the measure $\nu$ is normalized in such a way that generalized Weil's formula
$$
\int_Gg(x)dx =\int_{G/K}\left(\int_Kg(xk)d\beta(k)\right)d\lambda(\dot x) \eqno(1)
$$
holds for all $g\in L^1(G)$, where $\lambda$ denotes some left $G$-invariant measure on
$G/K$ (see \cite[Chapter VII, \S 2, no. 5, Theorem 2 ]{Bourb} and especially remark c)
after this theorem).
\footnote{$G$-left invariance of $\lambda$ means that $\lambda(xE) = \lambda(E)$ for every Borel subset $E$ of $G/K$ and
for every $x\in G$. This measure is unique up to a constant multiplier.}
Here and below we write $dx$ instead of $d\nu(x)$ and $dk$
instead of $d\beta(k)$.   We shall write also $d\dot x$ instead of $d\lambda(\dot x)$.

The function $g : G \to \mathbb{C}$ is called \textit{right $K$-invariant} if $g(xk) = g(x)$ for
all $x\in G$, $k\in K$. For such a function we put $\dot g(\dot x) := g(x)$. This definition
is correct and, since  $\int_K dk = 1$, for $g\in L^1(G)$ formula (1) implies that
$$
\int_Gg(x)dx =\int_{G/K}\dot g(\dot x)d\dot x.  \eqno(2)
$$

The map $g\mapsto \dot g$ is a bijection between the set of all right $K$-invariant
functions on $G$ (all right $K$-invariant functions from $L^p(G)$) and the set of
all functions on $G/K$ (respectively functions from $L^p(G/K)=L^p(G/K, \lambda)$).

Let an automorphism $A\in \mathrm{Aut}(G)$ maps $K$ onto itself. Since
$$
A(\dot x) := A(xK) = \{A(x)A(k) : k\in K\} = A(x)K = \pi_K(A(x)),
$$
we get a homeomorphism $\dot A: G/K \to G/K$, $\dot A(\dot x) := \pi_K(A(x))$. Then for
every right $K$-invariant function $g$ on $G$ we have $\dot g(\dot A(\dot x)) = g(A(x))$.

We put
$$
\mathrm{Aut}_K(G) := \{\dot A : A\in  \mathrm{Aut}(G), A(K) = K\}.
$$

\textbf{Definition 1.} \cite{JBSU} Let $(\Omega,\mu)$ be a measure space, $(\dot A(u))_{u\in
\Omega}\subset \mathrm{Aut}_K(G)$
a family of homeomorphisms of $G/K$, and $\Phi\in L^1_{loc}(\Omega,\mu)$. For a function $f$
on $G/K$ we define a \textit{Hausdorff operator} on $G/K$ as follows
$$
(\mathcal{H}_{\Phi, \dot A}f)(\dot x) :=\int_\Omega \Phi(u)f(\dot A(u)(\dot x))d\mu(u).
$$

It was shown in \cite[Theorem 1]{JBSU}  that under the conditions of Definition 1 for $p\in[1,\infty]$
the following inequality holds
$$
\|\mathcal{H}_{\Phi, \dot A}\|_{L^p(G/K)\to L^p(G/K)}\le \int_{\Omega}|\Phi(u)|(\mathrm{mod} A(u))^{-1/p}d\mu(u).
$$

Our goal is to give conditions for boundedness of
 operators $\mathcal{H}_{\Phi, \dot A}$  on atomic Hardy spaces $H^1$ over $G/K$.

Following \cite{CMM} we shall assume  that the group
$G$  possesses the following
properties:

(LDP) (\textit{local doubling property}): for every $b\in \mathbb{R}_+$ there exists a constant
$D_b$ such that for every ball $B$ in $G$ with radius $r_B < b$ the following inequality
holds
$$
\nu(2B)\le D_b\nu(B);
$$
where $2B$ denotes the ball in $G$ with the same center and radius $2r_B$.

(AMP) (\textit{approximate midpoint property}): there exist $R_0 \in [0, 1)$ and $\beta\in(1/2,1)$
 such that for every pair of points $x, y \in G$ with $\rho(x, y)> R_0$ there
exists a ball $B$ containing $x$ and $y$ with radius $r_B <\beta\rho (x,y)$.

By \cite[Remark 2.3]{CMM} the property (LDP) implies that for each
$r\ge 2$ and for each $b\ge 0$ there exists a constant $C>0$ such that
$$
\nu(B')\le C\nu(B) \eqno(LD)
$$
for each pair of balls $B$ and $B'$, with $B\subset B'$, $r_B\le b$, and $r_{B'}\le\tau r_B$. In the
following $D_{\tau,b}$ denotes the smallest constant for which (LD) holds.

\textbf{Remark 1}. As was noted  in \cite[p. 2]{CMM} important examples of metric measure
spaces which are locally doubling (but not doubling)  are complete Riemannian manifolds with Riemannian distance $\rho$ and
Riemannian density $\nu$, and with Ricci curvature bounded from below, a class
which includes all Riemannian symmetric spaces of the noncompact type and
Damek--Ricci spaces. A was mentioned in the introduction, a connected non compact semisimple Lie group with finite center also possesses the property (LDP).
 It is well known also that   every complete metric   space with path metric possesses the approximate midpoint property.

Under  conditions (LDP) and (AMP) we prove the theorem on boundedness of Hausdorff operators on atomic Hardy spaces $H^1$ over homogeneous spaces of locally
compact groups.

First recall that a function $a$ on $G$ is an ($(1,\infty)$-)\textit{atom} if

(i) the support of $a$ is contained in a ball $B(x, r)$;

(ii) $\|a\|_\infty \le \frac{1}{\nu(B(x,r))}$;

(iii) $\int_G a(x)d\nu(x) = 0$.

In case $\nu(G)<1$ we shall assume $\nu(G)=1$. Then the constant function
having value $1$ is also considered to be an atom.

According to \cite{CMM} \textit{an} $H^1_b$  \textit{atom} is an atom supported in a ball of radius at most $b$. Using $H^1_b$ atoms for each $b> 0$ the spaces $H^1_b = H^1_b(G):= H^{1,\infty}_b(G)$
on the group $G$ are defined in \cite{CMM} in the same manner as in the case of spaces
of homogeneous type considered in \cite{CW}, the only difference being that it is
required that the balls involved have at most radius $b$.
Furthermore, due to \cite[Proposition 4.3]{CMM} for $b>R_0/(1-\beta)$ we have
$H^1_c = H^1_b$ for all $c>b$.  So, we put $H^1(G) :=H^{1,\infty}_b(G)$ for such $b$.\footnote{Real Hardy spaces over compact connected (not necessary quasi-metric) Abelian
groups were defined in \cite{Indag}.}
 In the
following the constant  $b > R_0/(1-\beta)$ will be fixed.

\textbf{Definition 2.} (cf. \cite{JBSU}). We define the \textit{Hardy space} $H^1(G/K)$ as a space
of such functions $f =\dot g$ on $G/K$ that $g$ admits an atomic decomposition of
the form
$$
g =\sum_{j=1}^\infty \alpha_ja_j,  \eqno(3)
$$
where $a_j$ are right $K$-invariant $H^1_b$ atoms,
and $\sum_{j=1}^\infty |\alpha_j| < \infty$. In this case,
$$
\|f\|_{H^1(G/K)} := \inf \sum_{j=1}^\infty |\alpha_j|,
$$
and infimum is taken over all decompositions above of $g$.

Thus a function $f = \dot g$ from $H^1(G/K)$ admits an atomic decomposition
$f =\sum_{j=1}^\infty \alpha_j \dot a_j$ such that
$\sum_{j=1}^\infty |\alpha_j| < \infty$, and  $\|f\|_{H^1(G/K)}= \|g\|_{H^1(G)}$.

\textbf{Proposition 1.} \textit{The space $H^1(G/K)$ is Banach.  If for some  $h\in H^1(G)$
the inequality    $\int_K h(k)dk\ne 0$ holds the space  $H^1(G/K)$ is nontrivial.}

The proof of this proposition is  similar to the proof of Proposition 2 in \cite{JBSU}.

We need the following two lemmas to prove our main result.

\textbf{Lemma 1.} \cite{AddJMAA}  \textit{Every automorphism $A\in \mathrm{Aut}(G)$ of a locally
compact
metrizable group $(G,\rho)$ is Lipschitz. Moreover, one can choose the Lipschitz
constant to be}
$$
L_A = \kappa_\rho \mathrm{mod} A,
$$
\textit{where the constant $\kappa_\rho$ depends on the metric $\rho$ only}.

%\textbf{Remark 1}. The proof of the previous lemma in  \cite{AddJMAA} shows that $\kappa_{\rho_1}=1$ for some
%left invariant metric $\rho_1$ on $G$ which is compatible with the topology of $G$.

\textbf{Lemma 2.} \cite{JMAA} \textit{Let $(X;m)$ be a measure space and
$\mathcal{F}(X)$
be some Banach space of $m$-measurable functions on $X$. Assume that the
convergence of a sequence strongly in $\mathcal{F}(X)$ yields the convergence of some
subsequence to the same function for $m$-almost all $x\in X$. Let $(\Omega,q,\mu)$ be $\sigma$-compact quasi-metric space with quasi-metric $q$ and positive Radon measure $\mu$,and $F(u, x)$ be
a function such that $F(u, \cdot) \in \mathcal{F}(X)$ for $\mu$-a.e. $u\in \Omega$
and the map  $u\mapsto F(u, \cdot):\Omega\to \mathcal{F}(X)$ is Bochner integrable with respect to $\mu$. Then for $m$-a.e. $x\in  X$}
$$
\left((B)\int_\Omega F(u, \cdot)d\mu(u)\right)(x)=
\int_\Omega F(u, x)d\mu(u).
$$

Now we are in position to state and prove the following

\textbf{Theorem 1.} \textit{Let $\Omega$
 be $\sigma$-compact quasi-metric space with positive Radon
measure $\mu$. Let $G$ be a locally compact group with left Haar measure $\nu$ and
$\rho$  a left invariant metric which is compatible with the topology of $G$ and
conditions (LDP) and (AMP) hold. Let for some constant $C_1>0$ the family
$(\dot A(u))_{u\in \Omega}\subset \mathrm{Aut}_K(G)$
 satisfies $\mathrm{mod} A(u) \ge C_1$. Then for $\Phi\in L^1(\Omega, \mu)$ the
Hausdorff operator $H_{\Phi, \dot A}$  is bounded on the real Hardy space $H^1(G/K)$ and
for  $\tau = \max(2, \kappa_\rho/C_1)$ and some  constant $\gamma_{\tau,b}>0$ (depending only on $b$ and $\tau$)
the next estimate is valid}
$$
\|\mathcal{H}_{\Phi, \dot A}\|\le \gamma_{\tau,b}\|\Phi\|_{L^1(\Omega,\mu)}.
$$

Proof. If we set $X = G/K$ and $m = \lambda$ the pair $(X,m)$ satisfies the
conditions of Lemma 1 with $H^1(G/K)$ in place of $\mathcal{F}(X)$. Indeed, let $f_n =
\dot g_n \in H^1(G/K)$, $f =\dot g\in H^1(G/K)$, and $\|f_n- f\|_{H^1(G/K)}\to 0 (n\to \infty)$.
Since
$$
\|f_n- f\|_{L^1(G/K)} =\int_{G/K}|\pi_K(g_n-g)|d\lambda
$$
$$
=\int_G|g_n(x)-g(x)|dx \le\|g_n-g\|_{H^1(G)} =\|f_n - f\|_{H^1(G/K)}\to  0
$$
there is a subsequence of $f_n$ that converges to $f$ $\lambda$-a.e.

Then Definition 2 and lemma 1 imply for $f\in H^1(G/K)$ that
$$
\mathcal{H}_{\Phi, \dot A}f=\int_\Omega \Phi(u)f\circ\dot A(u)d\mu(u)
$$
(the Bochner integral).

Therefore (below $f =\dot g$, $g \in H^1(G) = H^1_b (G)$)
$$
\|\mathcal{H}_{\Phi, \dot A}f\|_{H^1(G/K)}\le\int_\Omega |\Phi(u)|\|f\circ\dot A(u)\|_{H^1(G/K)}d\mu(u)
$$
$$
=\int_\Omega |\Phi(u)|\|g\circ A(u)\|_{H^1(G)}d\mu(u).\eqno(4)
$$
If $g$ has representation (3) then
$$
g\circ A(u) =\sum_{j=1}^\infty \alpha_ja_j\circ A(u).\eqno(5)
$$

We claim that $a'_{j,u} := (1/D_{\tau,b}) a_j\circ A(u)$ is an atom from $H^1_{\tau b}(G)$.
Indeed,
if $a_j$ is supported in $B(x_j, r_j)$ ($r_j < b$) lemma 1 implies that $a_j\circ A(u)$ is
supported in
$$
A(u)^{-1}(B(x_j, r_j)) \subset B\left(A(u)^{-1}x_j,\frac{\kappa_\rho}{\mathrm{mod}A(u)}r_j\right)\subset B(A(u)^{-1}x_j, \tau r_j).
$$
And thus $a'_{j,u}$ is supported in $B(A(u)^{-1}x_j, \tau r_j)$ ($\tau r_j< \tau b$).

Next, note that $\|a_j\|\le 1/\nu(B(x_j, r_j))$. Since $\nu$ and $\rho$ are left invariant,
condition (LD) yields that
$$
\nu(B(A(u)^{-1}x_j, \tau r_j) = \nu(B(x_j, \tau r_j)) \le D_{\tau,b}\nu(B(x_j, r_j)).
$$
It follows that
$$
\|a'_{j,u}\|_\infty =
\frac{1}{D_{\tau,b}}\|a_j\|_\infty\le
\frac{1}{D_{\tau,b}(B(x_j, r_j))}\le
\frac{1}{\nu(B(A(u)^{-1}x_j,\tau r_j))}
$$
and $a'_{j,u}$ satisfies (i) and (ii) from the definition of an atom. The property
(iii) follows from the equality
$$
\int_Ga'_{j,u}d\nu = (1/D_{\tau,b})\mathrm{mod}A(u)^{-1}\int_Ga_jd\nu.
$$
We conclude that the function
$$
g \circ A(u) =\sum_{j=1}^\infty (\alpha_jD_{\tau,b})a'_{j,u}
$$
belongs to $H^1_{\tau b}(G)$ and
$$
\|g \circ A(u)\|_{H^1_{\tau b}}\le D_{\tau,b} \sum_{j=1}^\infty |\alpha_j|.
$$
So $\|g \circ A(u)\|_{H^1_{\tau b}}\le D_{\tau,b}\|g\|_{H^1_b}$.

On the other hand, by \cite[Proposition 4.3]{CMM} for $b > R_0/(1-\beta)$ we get
$H^1_{\tau b} = H^1_b$ and for some constant $C_{\tau,b} > 0$ depending only on $b$ and $\tau$
$$
\|g \circ A(u)\|_{H^1_{\tau b}}\le\|g \circ A(u)\|_{H^1_{b}}\le C_{\tau,b}\|g \circ A(u)\|_{H^1_{\tau b}}.
$$
Then
$$
\|g \circ A(u)\|_{H^1_{b}}\le C_{\tau,b}
\|g \circ A(u)\|_{H^1_{\tau b}}\le C_{\tau,b}D_{\tau,b}\|g\|_{H^1_{b}}.
$$
and thus
$$
\|\mathcal{H}_{\Phi,\dot A}f\|_{H^1(G/K)}=\|\mathcal{H}_{\Phi,\dot A}\dot g\|_{H^1(G/K)}\le\int_\Omega |\Phi(u)|\|\dot g\circ\dot A(u)\|_{H^1(G/K)}d\mu(u)
$$
$$
=\int_\Omega |\Phi(u)|\|g\circ A(u)\|_{H^1(G)}d\mu(u)\le C_{\tau,b}D_{\tau,b}\|\Phi\|_{L^1(\Omega,\mu)}\|g\|_{H^1(G)}
$$
$$
=C_{\tau,b}D_{\tau,b}\|\Phi\|_{L^1(\Omega,\mu)}\|f\|_{H^1(G/K)}
$$
and the proof is complete.

In \cite{JBSU} \'{C}esaro  operator over a homogeneous space $G/K$ was defined in the following way:
$$
(\mathcal{C}_{\dot A,\mu}f)(\dot x):=\int_{\{\mathrm{mod}A(u)\ge 1\}}\frac{f(\dot A(u)(\dot x))}{\mathrm{mod}A(u)}d\mu(u).
$$

\textbf{Corollary 1.} \textit{Under the conditions of theorem 1 we have that}
$$
\|\mathcal{C}_{A,\mu}\|_{H^1\to H^1}\le \gamma_{\tau,b}\int_{\{\mathrm{mod}A(u)\ge 1\}}\frac{d\mu(u)}{\mathrm{mod}A(u)}.
$$

Indeed, this follows from theorem 1, since for the  \'{C}esaro  operator
$$
\Phi(u)=\frac{\chi_{\{\mathrm{mod}A(u)\ge 1\}}(u)}{\mathrm{mod}A(u)}
$$
($\chi_E$ denotes the indicator of the subset $E\subset \Omega$).

\section{Hausdorff operators over the hyperbolic plane}

In this section we give an answer to the question, posed in \cite{JBSU}.

 Let $\mathbb{H}^2$ be the open upper half plane of the complex plane with the hyperbolic metric (the Poincar\'{e} model of the Lobachevsky plane). The  group $G: = SL(2)= SL(2,\mathbb{R})$ acts isometrically and transitively on $\mathbb{H}^2$ by the rule
$$
\left(\begin{array}{cc}
a &b\\
c &d
\end{array}\right)
\cdot z = \frac{az + b}{cz + d}.
$$
Since the stabilizer in $SL(2)$ of $i$ is the (maximal compact) subgroup $K:=SO(2)=SO(2,\mathbb{R})$ of $SL(2)$,
one can  identify $\mathbb{H}^2$ with the homogeneous space $SL(2)/SO(2)$ via the map $z = x\cdot i\mapsto \dot x:=\pi_K(x)=x SO(2)$ ($x\in SL(2)$) (this is a  diffeomorphism
of $\mathbb{H}^2$ onto $G/K$; see, e.g., \cite[Chapter III, \S 1]{Leng}). It is easy to verify that for $z\in\mathbb{H}^2$
$$
z=x(z)\cdot i
$$
where
 $$
 x(z)=\frac{1}{\sqrt{\mathrm{Im} z}}
\left(\begin{array}{cc}
\mathrm{Im} z &\mathrm{Re} z\\
0 &1
\end{array}\right).
$$

We  shall identify $z=x(z)\cdot i\in \mathbb{H}^2$ with $\pi_K(x(z))$.

It is known that $x\mapsto g^{-1}xg$ with $g\in GL(2,\mathbb{R})$ is the general form of automorphisms of the group $SL(2,\mathbb{R})$. Next, since $K=SO(2)$ is the connected component of the unit in $O(2)$, $K$ is a normal subgroup of $O(2)$. In other wards every automorphism $A(u)(x):=u^{-1}xu$ ($u\in O(2)$) of $SL(2)$ maps $K$ onto itself. By the previous  identification, $\dot A(u)(\dot x(z))=\pi_K(A(u)(x(z)))=(u^{-1}x(z)u)\cdot i$.

Hence  for our $G$, $K$, and $\Omega =O(2)$ the Definition 1 takes the form (we put $x=x(z)$ in this definition and  identify $\dot x(z)$ with $z$)
$$
(\mathcal{H}_{\Phi}f)(z) :=\int_{O(2)} \Phi(u)f((u^{-1}x(z)u)\cdot i)d\mu(u) \eqno(6)
$$
where $\mu$ stands for a (regular Borel) measure on $O(2)$ and $f$ is a  function on $\mathbb{H}^2$.

 Since $\mathrm{mod}A(u)=\Delta_G(u)$ \cite[Chapter VII, \S 1, n. 4]{Bourb} and $SL(2)$ is unimodular \cite[Chapter VII, \S 3, n. 3]{Bourb}, we have $\mathrm{mod}A(u)=1$ for all $u$. Let $\Phi\in L^1(O(2),\mu)$. Then theorem 1 from \cite{JBSU} yields, that the operator  (6) is bounded in $L^p(\mathbb{H}^2)$ for $p\in[1,\infty]$ and $\|\mathcal{H}_{\Phi}\|_{L^p\to L^p}\le \|\Phi\|_{L^1}$.

  Remark 1 implies that the group $SL(2)$ endowed with the path (Riemannian)  metric satisfies the conditions (LDP) and (AMP). Since  $\mathrm{mod}A(u)=1$, theorem 1 can be applied to the group $SL(2)$ with $C_1=1$. So
if $\Phi\in L^1(O(2),\mu)$ the operator  (6) is bounded on $H^1(\mathbb{H}^2)$ and $\|\mathcal{H}_{\Phi}\|_{H^1\to H^1}$ $\le \gamma\|\Phi\|_{L^1(\mu)}$ for some absolute constant $\gamma=\gamma_{\tau,b}$ ($\tau=\max(2,\kappa_{\rho})$,  $b>0$ is sufficiently large).

 It is well known that any   matrix from $O(2)$ looks like
 $$
k(\theta)=
\left(\begin{array}{cc}
 \cos\theta& -\sin\theta\\
\sin\theta &\cos\theta
\end{array}\right),\quad \theta\in[0,2\pi)
$$
(this matrix presents a rotation by $\theta$ by the origin in the Euclidean plain),

 or like:
$$
v(\theta)=
\left(\begin{array}{cc}
 \cos\theta& \sin\theta\\
\sin\theta &-\cos\theta
\end{array}\right),\quad \theta\in[0,2\pi)
$$
(this matrix presents a reflection in the Euclidean plain  across a line at an angle of $\theta/2$).

Consider both of these cases.

1) If $u=k(\theta)$ formula (6) takes the form
$$
(\mathcal{H}_{\Phi_1}f)(z) =\int_{0}^{2\pi} \Phi_1(\theta)f((k(\theta)^{-1}x(z)k(\theta))\cdot i)d\mu_1(\theta)
$$
where $\mu_1$ stands for a regular Borel measure on $[0,2\pi)$.

But, since $k(\theta)\cdot i=i$, we have
$$
(k(\theta)^{-1}x(z)k(\theta))\cdot i=k(\theta)^{-1}\cdot(x(z)\cdot(k(\theta)\cdot i))=k(-\theta)\cdot z=\frac{z\cos\theta+\sin\theta}{-z\sin\theta+\cos\theta}.
$$
Note that the M\"{o}bius transformation in the right-hand side induces a hyperbolic rotation of the half-plane $\mathbb{H}^2$ by the angle $2\theta$ about $i$ (see, e.~g., \cite[Lemma 9.19]{St}). So if $\theta$ runs over $[0,2\pi)$, the point $\frac{z\cos\theta+\sin\theta}{-z\sin\theta+\cos\theta}$ runs twice over the hyperbolic circle centered at $i$,  which  passes through  $z$ (this line is an Euclidean circle, too, see, e.~g., \cite[Corollary 5.3]{St}).

So  in this case the Hausdorff operator looks as follows:
$$
(\mathcal{H}_{\Phi_1}f)(z) =\int_{0}^{2\pi} \Phi_1(\theta)f(k(-\theta)\cdot z)d\mu_1(\theta)
$$
$$
=\int_{0}^{2\pi} \Phi_1(\theta)f\left(\frac{z\cos\theta+\sin\theta}{-z\sin\theta+\cos\theta}\right)d\mu_1(\theta).\eqno(7)
$$

2) Let $u=v(\theta)$.
Since $v(\theta)\cdot i=-i$ and $v(\theta)^{-1}=v(\theta)$, we have
$$
(v(\theta)^{-1}x(z)v(\theta))\cdot i=v(\theta)\cdot(x(z)\cdot(v(\theta)\cdot i))=-v(\theta)\cdot z=\frac{z\cos\theta+\sin\theta}{-z\sin\theta+\cos\theta}
$$
and we arrived at the same expression  for a Hausdorff operator as in the previous  case.

Thus, formula (7) gives us the general form of a Hausdorff operator on the hyperbolic plane.
This operator is bounded on  $L^p(\mathbb{H}^2)$ and $H^1(\mathbb{H}^2)$ if $\Phi_1\in L^1([0,2\pi),\mu_1)$ and by \cite[Theorem 1]{JBSU} and theorem 1
$$
\|\mathcal{H}_{\Phi_1}\|_{L^p\to L^p}\le \|\Phi_1\|_{L^1(\mu_1)},\quad
\|\mathcal{H}_{\Phi_1}\|_{H^1\to H^1}\le \gamma\|\Phi_1\|_{L^1(\mu_1)}.\eqno(8)
$$

\textbf{Example 2.} Consider the \'{C}esaro  operator over the hyperbolic plane (see corollary 1).
In our case we have  $\mathrm{mod}A(u)=1$ for all $u$. So we put $\Phi_1(\theta)=1$ in formula (7)
and define the (generalized)  \textit{\'{C}esaro  operator over the hyperbolic plane} by the formula
$$
(\mathcal{C}_{\mu_1}f)(z):=\int_{0}^{2\pi} f\left(\frac{z\cos\theta+\sin\theta}{-z\sin\theta+\cos\theta}\right)d\mu_1(\theta)
$$
(perhaps the best choice here is $d\mu_1(\theta)=d\theta/2\pi$). It follows from (8) that this operator is bounded on  $L^p(\mathbb{H}^2)$ and $H^1(\mathbb{H}^2)$ if  the measure $\mu_1$ is finite and
$$
\|\mathcal{C}_{\mu_1}\|_{L^p\to L^p}=\|\mu_1\|,\quad
\|\mathcal{C}_{\mu_1}\|_{H^1\to H^1}\le \gamma\|\mu_1\|.
$$

\end{document}